\documentclass[a4paper,11pt]{article}

\usepackage[utf8]{inputenc}
\usepackage[T1]{fontenc}
\usepackage[english]{babel}

\usepackage{amsmath}
\usepackage{amssymb}
\usepackage{amsfonts}
\usepackage{amsthm}
\usepackage{mathtools}
\usepackage{stmaryrd}
\usepackage{dsfont}
\usepackage{textcomp}

\usepackage{graphicx}
\usepackage{geometry}
\usepackage{csquotes}
\usepackage{enumerate}
\usepackage{todonotes}

\usepackage[capitalise]{cleveref}

\newtheorem{lemma}{Lemma}
\newtheorem*{lemma*}{Lemma}
\newtheorem{theorem}{Theorem}

\theoremstyle{definition}
\newtheorem{remark}{Remark}
\newtheorem{definition}{Definition}

\newtheorem{corollary}{Corollary}

\theoremstyle{plain}

\newtheorem*{notation*}{Notation}

\crefname{proposition}{Proposition}{Propositions}
\crefname{theorem}{Theorem}{Theorems}
\crefname{corollary}{Corollary}{Corollaries}
\crefname{korollar}{Corollary}{Corollaries}
\crefname{remark}{Remark}{Remarks}
\crefname{lemma}{Lemma}{Lemmas}
\crefname{definition}{Definition}{Definitions}
\crefname{example}{Example}{Examples}
\crefname{equation}{Equation}{Equations}
\crefname{aufgabe}{Exercise}{Exercises}
\crefname{abschaetzung}{Inequality}{Inequalities}

\newcommand{\R}{\ensuremath{\mathbb{R}}}
\newcommand{\N}{\ensuremath{\mathbb{N}}}

\newcommand{\Z}{\ensuremath{\mathcal{Z}}}

\newcommand{\Halphal}{
    \ensuremath{\mathcal{H}^{\alpha}([0,L]^2)}
}
\newcommand{\Halpha}{\ensuremath{\mathcal{H}^{\alpha}}}
\newcommand{\Qcal}{\mathcal{Q}}
\newcommand{\eps}{\varepsilon}
\renewcommand{\epsilon}{\varepsilon}
\newcommand{\T}{\mathbb{T}}
\newcommand{\ZZ}{\mathcal{Z} }
\newcommand{\cZ}{\ZZ}
\newcommand{\dd}{\,\mathrm{d}}

\def\Xint#1{%
    \mathchoice
        {\XXint\displaystyle\textstyle{#1}}%
        {\XXint\textstyle\scriptstyle{#1}}%
        {\XXint\scriptstyle\scriptscriptstyle{#1}}%
        {\XXint\scriptscriptstyle\scriptscriptstyle{#1}}%
    \!\int
}
\def\XXint#1#2#3{%
    {\setbox0=\hbox{$#1{#2#3}{\int}$}%
    \vcenter{\hbox{$#2#3$}}\kern-.6\wd0}%
}

\def\dashint{\Xint-}

\author{Dirk Bl\"omker, Johannes Rimmele}
\title{Stabilization by Rough Noise for an Epitaxial Growth Model}
\date{\today}

\begin{document}

\maketitle

\begin{abstract}
In this article, we study a model of epitaxial thin-film growth.
It was originally introduced as a phenomenological model of growth in the presence of a Schwoebel barrier, where diffusing particles on a terrace are not allowed to jump down at the boundary.

Nevertheless, we show that the presence of arbitrarily small space-time white noise due to fluctuations in the incoming particles surprisingly eliminates all nonlinear interactions in the model and thus has the potential to stabilize the dynamics and suppress the growth of hills in these models. 
\end{abstract}

\section{Introduction}
\label{sec:1}

Our starting point is an epitaxial thin-film  model for the growth of 
a crystalline surface in the presence of a Schwoebel barrier, which does not allow deposited atoms diffusing on the surface to jump down steps.

For a graph $h(t,x)$ of a surface at time $t>0$ over a point $x$, 
the model is given by the following stochastic partial differential equation (SPDE) 
\begin{equation}
\label{e:HS}
\partial_t h = - \delta \Delta^2 h - \nabla \cdot \frac{\nabla h}{1+|\nabla h|^2} + \sigma \xi 
\end{equation}
where $\xi=\partial_t W$ is space-time white noise given by the derivative of a standard cylindrical Wiener process $W$. The noise arises 
naturally due to thermal fluctuations in the incoming particles.  

Moreover, $\delta$ and $\sigma$ are small positive parameters.
In the physics literature  \eqref{e:HS} is usually considered either on the whole space $\mathbb{R}^d$, 
with dimension $d\in\{1,2\}$ or with periodic boundary conditions on the torus $\mathbb{T}^d$, where the interesting physical dimension is $d=2$.
The equation is posed for a growing surface in a moving frame, thus we assume in the bounded domain case 
$\int_{\mathbb{T}^d} h(t,x)dx=0$ for all $t\geq0$.

Originally, the model was introduced by Hunt, Sander, et al.\ 
for $d=1$ in \cite{HOWOS:94} and in $d=2$ in \cite{JOHGSSO:94}.
For general surveys of surface growth, see, for example, \cite{KS:95,BS:95} or \cite{KS:91,LDS:91},
where this equation is treated.

The equation is derived phenomenologically for the 
epitaxial growth of a crystalline surface where deposited atoms are allowed to diffuse on the terraces of the surface. 
The linear operator, i.e.\ the bilaplacian, arises from linearized surface diffusion,
while  the nonlinearity should have the following important features due to the Schwoebel barrier
that prohibits the diffusing atoms on the surface to jump down at the boundary of terraces,
which leads to an effective uphill current.
\begin{itemize}
	\item {\em Linear instability}: For a small gradient $\nabla h$ of the surface, the nonlinearity behaves like $-\Delta h$, which should lead to a linear instability and the growth of hills on the initially flat surface.  
	\item The nonlinearity vanishes if $|\nabla h|$ is large, as atoms immediately get stuck on a steep surface. 
\end{itemize}

There are numerous other models for these types of growth processes
in the physics 
literature where the nonlinearity is 
of the type $\nabla \cdot f(\nabla h)$ 
with $f:\mathbb{R}^d\to\mathbb{R}^d$.
See \cite{PGMPV:00} for examples. 
A model where for some  $\gamma>0$ one has 
\[|f(\nabla h)| \sim |\nabla h|^{-\gamma}
\] 
for large gradients can also be found in \cite{PT:00}. 
Other models in the literature consider functions $f$ that depend non-isotropically on $\nabla h$. This would allow for slope selection of the hills, see \cite{ZZ:13}, for an example. Here we focus on the simpler example \eqref{e:HS} on the torus.

\paragraph{Mathematical problem:} 
The key problem of the equation in the physically relevant case $d=2$ is the regularity of the solution. One expects that with space-time white noise
the solution is not regular enough to make sense out of the gradient $\nabla h$ as a function.
From the linearization one expects the surface $h$ to be at most H\"older 
with exponent strictly less than one.

\paragraph{Deterministic surface growth model:}
While the stochastic equation does not appear to have been studied in the mathematical literature,
the deterministic equation (i.e., \eqref{e:HS} with $\sigma=0$)
and various modifications have been studied in some detail in recent years.

Well-posedness and spectral Galerkin methods 
were studied in \cite{LL:03,LL:05,Li:06},
where the dynamics of phase separation was also investigated both 
numerically and through asymptotic scaling.
For results on existence and uniqueness in a more general model, 
see, for example, \cite{Ag:15}.
For our model, due to the global Lipschitz continuity of the nonlinearity as a map, for example, from
the Sobolev space $H^{1}(\mathbb{T}^d)$ into $H^{-1}(\mathbb{T}^d)$ 
there are no serious problems for the well-posedness of the equation. 

The error estimates for various numerical discretization methods were studied in numerous publications.
We mention \cite{LSL:17,CCWWW:12,QSZ:15,ZZ:13} for a few non-exhaustive examples, where also more general variants of \eqref{e:HS} were studied. 

The long-time behavior of the deterministic version of 
\eqref{e:HS} was studied 
in several publications \cite{DX:19,GMY:11,AMY:17}. 
The solution of the initial value problem for the 
deterministic equation generates a dissipative
dynamical system that exhibits an attractor
dominating the long-time behavior. Moreover,
the existence of an energy 
\[
E(u)= \frac12 \int_{\mathbb{T}^d}  \delta |\Delta h|^2 - \ln(1+|\nabla h|^2) dx 
\]
shows that our model is a gradient flow.
Thus, it is well known that on a bounded domain 
the attractor consists only of equilibria and heteroclinic connections.
Note that due to Poincar\'e's inequality the energy is bounded from below,
but for small $\delta>0$, the flat surface, i.e.\ the constant $h=0$, is no longer the minimizer of the energy, 
and thus unstable, while there is an (up to translation) unique global minimizer.

However, the detailed structure of this attractor and the energy landscape
are far from being fully understood.
This also leads to interesting questions for the SPDE 
about the stochastic motion of hills, which have not yet been studied.  

\paragraph{Renormalization:}
The well-posedness theory for irregular SPDEs driven by excessively rough noise is by now well developed, ranging from the first ideas of Da Prato and Debussche \cite{DPD02,DPD03} 
to the theory of paracontrolled distributions by Perkowski, Gubinelli et al.\ \cite{GIP:15} and the celebrated theory of regularity structures \cite{Ha:14a,Ha:14b} (see also \cite{Ha:13}) by Martin Hairer. 

There are now numerous results in the literature on how to give 
meaning to equations that are ill-posed for direct methods like fixed-point arguments or Galerkin methods
due to the presence of spatial noise that is too rough,
and we do not try to give a complete overview here. 
But many of these results focus on equations with  polynomial nonlinearities,
and equation \eqref{e:HS} does not seem to fit into this setting.

Classical examples with quadratic nonlinearities mentioned above include Burgers-type, 
KPZ, and Navier--Stokes equations.   
There are also results for the stochastic quantization equation \cite{CC:18},
FitzHugh-Nagumo \cite{BK16}, the $\Phi^4_3$-model \cite{Ha:16},
or the dynamical sine--Gordon model \cite{HS:16}, where the nonlinearity is a sine function, 
but there are numerous other examples in the literature.

In a simplified setting, the basic idea of all these approaches for irregular SPDEs, which we will also utilize here, 
is a regularization of the noise which then gives a meaning to the equation in the 
limit of vanishing regularization. 
In this limit, diverging counterterms often 
have to be taken into account.

In our setting the regularity of the noise 
is only slightly too irregular to give a meaning of the equation directly via fixed-point arguments 
or Galerkin methods.
Thus, we shall rely on a straightforward Fourier-series computation. The underlying idea can be traced back to Da Prato and Debussche \cite{DPD02,DPD03}.

In the context of surface growth equations, this approach has already been used in \cite{BR:13}, where the equation has a different structure than \eqref{e:HS}. 
While the linearization around $0$ is the same, 
the
nonlinearity is quadratic. For a review of this equation, see \cite{BR:15}.

\paragraph{Stabilization - Triviality of SPDEs:}
Stabilization of stochastic differential equations 
or SPDEs via noise is a broad area of research 
and we are not trying to give an overview here.

The stabilization effect studied in this work is a mechanism in which stabilization arises 
from the spatial roughness of the noise. It is related to the work of Hairer, Ryser, and Weber 
on triviality for SPDEs \cite{HRW:12}.

In \cite{HRW:12}, the authors studied the Allen--Cahn equation, 
also called $\Phi^4$-model or Ginzburg-Landau equation, 
in two spatial dimensions with periodic boundary 
conditions. Here it is not possible via standard methods
to consider the equation 
with space-time white noise, as the solution is not regular enough to 
make sense out of the cubic nonlinearity.
One therefore usually considers regularized noise and 
renormalized nonlinearities 
with diverging constants.

Without renormalizing the nonlinearity, the regularized equation in \cite{HRW:12} is 
\begin{equation}\label{eq:HRW}
\partial_t u_\eps 
=  \Delta u_\eps + f(u_\eps) 
+  \sigma_\eps \Qcal_\eps^{1/2} \partial_t W,
\end{equation}
where the covariance operator $ \sigma_\eps \Qcal_\eps^{1/2} $ is split explicitly into the scalar noise strength~$\sigma_\eps$ and spatial correlation given by the symmetric operator $\Qcal_\eps$, which may, for example, be introduced through a Fourier cut-off of the noise. 

The nonlinearity is 
\[
f(u)=u-u^3
\]
consisting of a stabilizing cubic $-u^3$ and a
linearly unstable term $u$, which leads on large domains to 
a linear instability and thus to pattern formation 
starting from small solutions.

The key problem here in spatial dimension two is that 
for $Q_\eps=I$, i.e.\ for space-time white noise,
the solution is only a distribution and 
not regular enough 
to define the cubic or any other power in the nonlinearity $f$.   

In their work \cite{HRW:12}
the authors identified three regimes. 
If the noise strength~$\sigma_\eps$ is too small, the noise vanishes, while in an intermediate regime 
one obtains a non-trivial limit, 
where the limit $u$ of $u_\eps$ solves 
a deterministic Allen-Cahn equation with an additional stabilizing 
term $-Cu$ for a constant $C>0$. 

The case of particular interest here 
is that of fixed noise strength $\sigma_\eps=\sigma$,
where, through averaging, the interaction of the noise with the cubic term
generates a divergent 
stabilizing linear term $-C_\eps u_\eps$ 
with $C_\eps\to \infty$ for $\eps\to 0$.

Thus, the solution $u_\eps$ would converge to zero 
removing  all possible  instabilities present in $f$.
But the convergence is only in a Sobolev or Besov space of negative order in derivatives, 
as even for small but positive $\eps>0$ the solution $u_\eps$ is highly oscillatory in space due 
to the noise. Thus, the roughness of the noise basically 
dominates all deterministic dynamics present 
in the model. This differs from our equation, in which only the gradient of the solution is ill-defined.     

A similar result in the same spirit was obtained recently 
by Oh, Okamoto and Robert in \cite{OOR:20} 
for the stochastic nonlinear wave equation.

In \cite{BT22} 
we studied results of this type 
not only for the Allen-Cahn equation 
but in a more general setting, 
where the linear operator 
and the nonlinearity $f$ may also change with $\eps$. 
An interesting result is \cite{HW:15} by Hairer and Weber
where the authors again study, for the Allen--Cahn equation, the limits of the noise 
strength tending to zero and the loss of regularity 
of the noise both separately and combined. 

In Section \ref{sec:2}, we state the setting and the main results.
The properties of the stochastic convolution $Z$ and its regularization $Z_\varepsilon$ are covered in Section \ref{sec:3}, including a convergence result, while Section \ref{sec:4} provides a standard existence and uniqueness result for solutions $u_\varepsilon$ of the regularized problem.
A crucial result is the deterministic uniform bound for the difference $u_\eps-Z_\eps$
in Section \ref{sec:5}. This allows us to prove the vanishing of the nonlinearity in Section \ref{sec:6} and the main result in Section \ref{sec:7}.


\section{Setting and Main Results}
\label{sec:2}

This article is devoted to the study of a random solution~$ u : [0,T] \times [0,L]^2 \to \R $
of the stochastic partial differential equation
\begin{align} \label{SPDE}
\partial_t u(t,x) = - \Delta ^2 u(t,x)  - \nabla  \cdot  \frac{\nabla u(t,x)}{1 + \left | \nabla u(t,x) \right | ^2} \
+ \sigma \ \partial _t  W(t,x)
\end{align}
for~$x \in [0,L]^2$ and $t \in [0,T]$. 
We assume that the solution is expressed in a moving frame such that the mean-value condition
\begin{align*}
\dashint_{[0,L]^2} u(t,x) dx = 0
\end{align*}
holds for each $t \in [0,T]$.
For notational simplicity, we set $\sigma=\delta=1$ throughout the analysis.
All arguments extend to arbitrary fixed $\sigma,\delta>0$, with the constants depending on $\sigma$ and $\delta$.
Furthermore, we assume that $u$ satisfies periodic boundary conditions, which yield the following explicit Fourier-series expansion:
$$
u(t, x) =
\sum_{k \in \mathbb{Z}^2} u_k(t) \ e_k(x),
$$
for $x = \left( x_1, x_2 \right) \in [0,L]^2$,
where $\left(e_k \right)_{k \in \mathbb{Z}^2 }$
is an orthonormal basis of $L^2 \left([0,L]^2 \right)$
given by~$
e_k(x)
\coloneqq
\omega_{k_1}(x_1) \ \omega_{k_2}(x_2)
$
and
\begin{align*}
\omega_{k_i} (z) \coloneqq
\begin{cases}
C \sin \left(\frac{2 \pi k_i}{L} z  \right),& \quad  \text{ for $k_i >0$,}\\
\frac{1}{\sqrt{L}}  \quad \quad \quad ,& \quad \text{ for $k_i = 0$,}\\
C \cos \left(\frac{2 \pi k_i }{L} z\right),& \quad \text{ for $k_i <0$}
\end{cases}
\end{align*}
for $i \in \{1,2\}$ and $C \coloneqq	\sqrt{\frac{2}{L}} $.
Using these eigenfunctions and the identity
\begin{align*}
- \Delta e_k = \left( \frac{2 \pi}{L }\right)^2 |k|^2 \ e_k
\end{align*}
for each $k \in \mathbb{Z}^2 \setminus \{0\}$. We define the eigenvalues~$\mu_k \coloneqq \left(\frac{2 \pi |k|}{L} \right)^2$ of the Laplace operator and derive the eigenvalues of the \emph{Bilaplace operator}, 
according to
\begin{align*}
 \Delta^2 e_k =   \left( \frac{2 \pi}{L }\right)^4 |k|^4 \ e_k
= 
 \mu_k ^2 \ e_k.
\end{align*}
For simplicity of notation, we define
$\mathcal{Z} \coloneqq \mathbb{Z}^2 \setminus \{0\}$ and characterize the \emph{fractional Sobolev space} $\mathcal{H}^\alpha = W^{\alpha , 2}$ for every $\alpha \in \R$ by
\[
\Halpha \coloneqq
\mathcal{H}^\alpha \left([0,L]^2 \right)
\coloneqq 
\left \{
u(x) = 
\sum_{k \in \mathcal{Z}}
u_k \ e_k \ 
\colon
\sum_{k \in \mathcal{Z}}
\mu_k ^\alpha \
u_k^2  < \infty 
\right \}
\]
and the \emph{fractional Laplace operator}
$(- \Delta )^{\alpha /2}
\colon
\Halphal \to L^2 \left([0,L]^2\right)
$
by
\begin{align*}
\sum_{k \in \mathcal{Z}}
u_k e_k
\mapsto
\sum_{k \in \mathcal{Z}}
\mu_k ^{\alpha/2} \
u_k e_k\;.
\end{align*}
We equip the space~$\mathcal{H}^\alpha \left([0,L]^2 \right)$ with the norm~$\| . \|_{\mathcal{H}^\alpha}$, defined for $u \in \mathcal{H}^\alpha \left([0,L]^2 \right)$ as follows:
\begin{align*}
\| u \|_{\mathcal{H}^\alpha}
\coloneqq
\left \| 
(- \Delta )^{\alpha /2} \ u 
\right\|_{L^2} =
\left(
\sum_{k \in \mathcal{Z}}
\mu_k ^\alpha \
u_k^2
\right)^{1/2}.
\end{align*}

\noindent 
For vector-valued functions, we set
\[
\mathcal{H}^{\alpha}(\T^2,\R^2)
\coloneqq
\mathcal{H}^{\alpha}(\T^2)\times\mathcal{H}^{\alpha}(\T^2),
\]
equipped with the norm
\[
\|(g_1 , g_2 )\|_{\mathcal{H}^{\alpha}(\T^2,\R^2)}
\coloneqq
\left(
\|g_1\|_{\mathcal{H}^{\alpha}(\T^2)}^2
+
\|g_2\|_{\mathcal{H}^{\alpha}(\T^2)}^2
\right)^{1/2}.
\]

\begin{remark}
By Poincaré's inequality, on the space of functions with zero spatial mean, the norm $\| \cdot \|_{\mathcal{H}^\alpha}$ 
is equivalent to the standard~$H^\alpha$-norm.
\end{remark}

\subsection{Analytic Semigroup}

As the operator $A \coloneqq -\Delta ^2$ is strongly elliptic on~$L^p \left([0,L]^2\right)$ for every $p \in (1,\infty)$, it follows that
$\left(e^{tA} \right)_{t \geq0}$ is an analytic semigroup on~$L^p \left([0,L]^2\right)$ for every $p \in (1,\infty)$ \cite[p. 212-214]{PAZY}.
Using the Fourier expression for $u \in L^2\left([0,L]^2\right)$, we obtain
\begin{align*}
e^{tA} u(x) = 
\sum_{k \in \mathcal{Z}}  e^{-t \mu_k^2} u_k e_k(x)
\end{align*}
for $t \geq 0$ and $x \in [0,L]^2$.
We now consider the operator norm of 
$$\left(e^{tA} \right)_{t \geq0} : \mathcal{H} ^{\alpha} \left([0,L]^2 \right) \to \mathcal{H} ^{\beta } \left([0,L]^2 \right)
$$
where $\beta > \alpha$.
For all $ u \in \mathcal{H}^{\alpha} \left([0,L]^2\right)$  we have the following well-known inequality for each $t >0$
\begin{align}
\label{e:SG}
\left \|
e^{tA}
\right \|_{L \left( \mathcal{H}^{\alpha},  \mathcal{H}^{\beta } \right)}
\leq \left ( \frac{\beta -\alpha}{ 4e}\right ) ^{\frac{\beta - \alpha }{4}} t^{\frac{\alpha -\beta}{4}}.
\end{align}
Later, integrability requires the order of the singularity in $t$ to be strictly less than $1$.
Due to our choice of the Sobolev norm in $\mathcal{H}^\alpha$, the operator
$$
\nabla : \Halphal \to \mathcal{H} ^{\alpha -1 } \left([0,L]^2, \mathbb{R}^2 \right)
$$
is an isometric bounded linear operator. Using the integration-by-parts formula,
for $s \in \mathcal{H}^\alpha$ we obtain
\[ \|\nabla s\|_{\mathcal{H}^{\alpha-1}
\left(
[0,L]^2, \R^2
\right)
}^2 = \langle  (-\Delta) s, s\rangle_{\mathcal{H}^{\alpha-1}} = \|(-\Delta)^{1/2} s\|_{\mathcal{H}^{\alpha-1}}^2
=  \| s\|_{\mathcal{H}^\alpha}^2.
\]
For the divergence, we obtain 
\[
\|\nabla \cdot g\|_{L^2} = \| \partial_x g_1+ \partial_y g_2\|_{L^2} \leq C \|g\|_{\mathcal{H}^1\left( [0,L]^2, \R^2 \right)},
 \]
where $g$ belongs to $\mathcal{H}^1\left( [0,L]^2, \R^2 \right)$.
Thus,
$$e^{tA}\nabla  \cdot  \colon
\mathcal{H}^{\alpha -1} \left([0,L]^2, \mathbb{R}^2  \right) \rightarrow \mathcal{H}^{\alpha} \left([0,L]^2\right)
$$
is a bounded linear operator with 
\[
\left \|
e^{tA} \nabla \cdot
\right \|_{L \left( \mathcal{H}^{\alpha-1} \left([0,L]^2 , \R^2 \right),  \mathcal{H}^\alpha \right)}
\leq 
\left \|
e^{tA}
\right \|_{L \left( \mathcal{H}^{\alpha-2},  \mathcal{H}^\alpha \right)}
\left \|
\nabla \cdot\right \|_{L \left( \mathcal{H}^{\alpha-1} \left([0,L]^2 , \R^2 \right),  \mathcal{H}^{\alpha-2} \right)}
\leq 
C \ t^{-1/2}.
\]

\subsection{White Noise}

In this subsection, we state the basic properties of \emph{space-time white noise}~$\partial_t W$ 
 and 
a cylindrical Wiener process
$$W(t,x) \coloneqq \sum_{k \in \mathcal{Z}} \beta_ k(t) e_k(x) ,
$$
for $x \in [0,L]^2$ and $t \in [0,T]$,
where $\left(\beta_k\right)_{k \in \mathcal{Z}}$ 
is a sequence of i.i.d.\ real-valued Wiener processes. We define the stochastic convolution
\begin{equation}
\label{e:defZ}
Z (t,x) \coloneqq 
\int_{0}^{t} e^{(t-s) A} dW (s,x) 
= \sum_{k \in \mathcal{Z} }
\int_{0}^{t} e^{-  (t-s) \mu_k^2} \ d\beta_k (s) \ e_k(x)
\end{equation}
which characterizes the solution of 
$
	\partial_t u  = - \Delta ^2 u + \partial_t W
$
with initial condition $u_0 = 0$.
We now analyze whether $\left(Z(t)\right)_{t \in [0,T] }$ belongs to $\mathcal{H}^\alpha
\left([0,L]^2 \right)
$.  It is well known that a short calculation using Fourier series shows
\[
\mathbb{E} \left[
\left \|
Z (t)
\right \|^2 _ {\mathcal{H}^1}
\right]
 = 
\infty  
\]
for each $t >0$,
but for~$\alpha < 1$ we get 
\[
\mathbb{E} \left[
\left \|
Z (t)
\right \|^2 _ {\mathcal{H}^\alpha}
\right]
= 
\sum_{k \in \mathcal{Z}} 
\frac{1- 
e^{-2t \mu_k ^2}}{2 \mu_k ^{2-\alpha}}
<
\infty
\]
for each $t \in [0,T]$.
For the analysis of the regularized equation \cref{e:SPDEreg}, we require
the regularity $\mathcal H^1([0,L]^2)
$
for every \(t\in[0,T]\). This regularity is needed in our fixed-point
framework to ensure that \(\nabla u_\varepsilon(t)\) and hence the
nonlinearity are well defined.

\subsection{Regularization of the stochastic solution}

Let us now define a regularized stochastic convolution.
\begin{definition}[Regularization coefficients]
\label{def:regcoeff}
Let
 $
\left(\alpha^{(\varepsilon)}\right)_{\varepsilon>0}$ 
with $ 
\alpha^{(\varepsilon)}
\coloneqq
\left(\alpha_k^{(\varepsilon)}\right)_{k\in\mathcal Z}$
be a family of real-valued sequences satisfying the following
assumptions:
\begin{enumerate}[(i)]
	\item The coefficients are uniformly bounded:
	$ \displaystyle
	\sup_{\varepsilon>0}
	\sup_{k\in\mathcal Z}
	\left|\alpha_k^{(\varepsilon)}\right|
	<\infty.
	$ 
	\item For every $\varepsilon>0$, there exist
	$\delta(\varepsilon)>0$ and $C_\varepsilon>0$ such that
	\[
	\left|\alpha_k^{(\varepsilon)}\right|
	\leq
	C_\varepsilon |k|^{-\delta(\varepsilon)},
	\qquad k\in\mathcal Z.
	\]
	\item For every fixed $k\in\mathcal Z$, \quad
	$ 
	\alpha_k^{(\varepsilon)}
	\xrightarrow{\varepsilon\to0}
	1.
	$
	\item The coefficients are radially symmetric, that is, for every
	$\varepsilon>0$ there exists a function
	$a_\varepsilon:[0,\infty)\to\mathbb R$ such that
	$ 
	\alpha_k^{(\varepsilon)}
	=
	a_\varepsilon(|k|)$ for $ k\in\mathcal Z.
	$ 
    \end{enumerate}
\end{definition}

\begin{definition}[Regularized noise]
\label{def:Weps}
Let the coefficients
$\left(\alpha_k^{(\varepsilon)}\right)_{k\in\mathcal Z}$
satisfy \Cref{def:regcoeff}. For every $\varepsilon>0$, define the
bounded linear operator $R_\varepsilon$ by
\[
R_\varepsilon e_k
\coloneqq
\alpha_k^{(\varepsilon)}e_k,
\qquad k\in\mathcal Z.
\]
The regularized cylindrical Wiener process is defined by
\begin{equation}
\label{eq:Weps}
W_\varepsilon(t)
\coloneqq
R_\varepsilon W(t)
=
\sum_{k\in\mathcal Z}
\alpha_k^{(\varepsilon)}
\beta_k(t)e_k,
\qquad t\geq0,
\end{equation}
where the series is understood in the cylindrical sense.
\end{definition}

\begin{definition}[Regularized stochastic convolution]
\label{def:Zeps}
For every $\varepsilon>0$, define the regularized stochastic
convolution by
\begin{equation}
\label{e:defZeps}
Z_\varepsilon(t)
\coloneqq
\int_0^t e^{(t-s)A}\,\dd W_\varepsilon(s),
\qquad t\in[0,T].
\end{equation}
Equivalently, for $t\in[0,T]$ and $x\in[0,L]^2$,
\[
Z_\varepsilon(t,x)
=
\sum_{k\in\mathcal Z}
\alpha_k^{(\varepsilon)}
\int_0^t
e^{-(t-s)\mu_k^2}\,\dd\beta_k(s)\,e_k(x).
\]
\end{definition}

Thus, using a straightforward Fourier-space calculation based on the Itô isometry, we derive
\[
\mathbb{E} \left[
\left \|
Z_\varepsilon (t)
\right \|^2 _ {L^2}
\right] 
= 
\sum_{k \in \mathcal{Z}} 
\left(\alpha_k^{(\varepsilon)}\right)^2
\left[
1- 
e^{-2t \mu_k ^2}
\right]
\frac{1}{2 \mu_k ^2 }
\leq 
\frac{1}{2}
\sum_{k \in \mathcal{Z}} \left(\alpha_k^{(\varepsilon)}\right)^2 \mu_k ^{-2} 
< \infty.
\]
Furthermore, we obtain
\begin{align*}
\mathbb{E} \left[
\left \|
\nabla
Z_\varepsilon (t)
\right \|^2 _ {L^2\left([0,L]^2, \mathbb{R}^2  \right)}
\right] 
=& \
\mathbb{E} \left[
\left \|
Z_\varepsilon (t)
\right \|^2 _ {\mathcal{H}^1}
\right] 
\\
= & \ 
\mathbb{E} \left[ \
\sum_{k \in \mathcal{Z}} \mu_k \ \left(\alpha_k^{(\varepsilon)}\right)^2 \left( \int_{0}^{t} e^{-(t-s)\mu_k ^2} d \beta_k (s)  \right)^2 
\right]
\\
=&
\sum_{k \in \mathcal{Z}} \left(\alpha_k^{(\varepsilon)}\right)^2 
\left[
1- 
e^{-2t \mu_k ^2}
\right]
\frac{1}{2 \mu_k }
\leq  C
\sum_{k \in \mathcal{Z}} \left(\alpha_k^{(\varepsilon)}\right)^2 
\frac{1}{ |k|^2 } 
 <  \ \infty.
\end{align*}
Thus, for all $\varepsilon>0$ the stochastic convolution $Z_\varepsilon$ has a well-defined spatial gradient.

\subsection{Main Results}

		\begin{definition}
			\label {def:mild}
			Let  $A \coloneqq -\Delta^2$ and 
			$
			\mathfrak{F} \left(u_\varepsilon  \right)
			\coloneqq - \nabla \cdot f(\nabla u_\varepsilon)
			\coloneqq - \nabla \cdot \frac{\nabla u_\varepsilon}{1 + \left | \nabla u_\varepsilon \right |^2}
			$
			 and $Z_\varepsilon$ be the stochastic convolution, as defined in \Cref{def:Zeps}.
			Let $u_0\in\mathcal{H}^{1  } (\T^2)$.
			For~$T>0$, we call a process $(u_\varepsilon(t))_{t\in[0,T]}$ with continuous paths in~$\mathcal{H}^1 (\T^2) $ such that for all $t\in[0,T]$ and $x \in \T^2$
			\begin{align*}
				u_\varepsilon (t,x) =
				e^{tA} u _0 (x) + 
				\int_{0}^t
				e^{(t-s)A} \mathfrak{F} \left(u_\varepsilon (s) \right)(x)\dd s+ Z_\varepsilon(t,x),
			\end{align*}
			a \emph{mild solution} of the regularized SPDE 
			\begin{equation}  \label{e:SPDEreg}
				\partial_t u_\varepsilon (t,x) = -\Delta^2 u_\varepsilon  (t,x)
				-
				\nabla \cdot \frac{\nabla u_\varepsilon (t,x)}{1 + \left | \nabla u_\varepsilon  (t,x)\right |^2}
				+
				\partial_t W_\varepsilon (t,x), \quad  u_{\varepsilon}(0)=u_0.
			\end{equation}
		\end{definition}
	
	\begin{theorem} \label {thm:conv}
        Let $u_0\in\mathcal{H}^{1  } (\T^2) \cap W^{1, \infty} (\T^2)$, $T>0$ and $u_\varepsilon$ be the mild solution of \eqref{e:SPDEreg}.
		Then, for each $p \geq 1$, $u_\varepsilon$ converges to $u$ if~$\varepsilon\to 0$ in~$L^p(\Omega,C^0([0,T]\times\T^2))$, where $u$ is a mild solution of
		\begin{equation}  \label {e:*}
			\partial_t u =-\Delta^2u + \partial_t W,
            \qquad
            u(0)=u_0.
		\end{equation} 
	\end{theorem}
	
	\noindent
	\begin{remark}
		Hill formation is expected for \eqref{e:SPDEreg} due to the linear instability indicated by the approximation
		$$- \nabla \cdot f \left(\nabla u \right) \approx - \Delta u$$
		 for small $\nabla u$.
		This approximation yields the linearized equation
		\begin{align*}
			\partial_t u = - \left( \Delta^2 + \Delta \right) u +   \sigma \partial_t W.
		\end{align*}
		For a more detailed analysis of the linear instability effects in the present model, see \cite{BR:GR26}; for a related discussion in the context of the Cahn--Hilliard equation, see \cite{DBMPTW:01}. In contrast, when the nonlinearity vanishes, no hills are observed in the solution of \eqref{e:*}, where solutions are expected to be merely of order $\sigma$.
		In particular,
		this yields \eqref{e:*}.
		However, we analysed this linearization approach in more detail in \cite{BR:GR26}.
	\end{remark}

\section{Stochastic Convolution}
\label{sec:3}

To apply Banach's fixed-point theorem, we first need the gradient of the stochastic convolution to be defined.
First, we recall that the Fourier series of~$Z(t,x)$ defined in \eqref{e:defZ} and~$Z_\varepsilon(t,x)$ defined in \eqref{e:defZeps} are well-defined pointwise as real-valued random variables for all $(t,x) \in [0,T] \times [0,L]^2$.

\begin{lemma} \label{lem:Zptwse}
For every~$x \in [0,L]^2$ and $t \in [0,T]$
the series defining $Z(t,x)$ and $Z_\epsilon(t,x)$ converge in $L^2(\Omega,\mathbb{R})$ and 
define well-defined real-valued Gaussian random variables.
\end{lemma}

\begin{proof}
We prove this statement only for $Z_\epsilon$ with $\epsilon > 0$,
as the proof for $Z$ is completely analogous.
 To show this, we consider the finite truncation 
\begin{align*}
Z_\epsilon^{(n)}(t,x) \coloneqq 	\sum_{k \in B_n(0)} \alpha_k^{(\varepsilon)} 
\int_{0}^{t} e^{-(t-\tau)\mu_k ^2} \ d \beta_k (\tau) \
e_k(x)
\end{align*}
by restricting the sum to the ball~$B_n(0) \coloneqq 
\left \{
y \in \mathcal{Z} :
|y| \leq n
\right \}
\subset \mathcal{Z}
$ with radius~$n \in \N^*$.
Now $Z^{(n)}_\varepsilon(t,x)$ is a well-defined real-valued Gaussian, and by the Itô isometry we obtain
\begin{align*}
\mathbb{E}
\left[
\left |
Z_\epsilon^{(n)}(t,x)
\right |^2
\right] 
= & \
\mathbb{E}
\left[
\left |
\sum_{k \in B_n(0)} \alpha_k^{(\varepsilon)} 
\int_{0}^{t} e^{-(t-\tau)\mu_k ^2} \ d \beta_k (\tau) \
e_k(x)
\right |^2
\right] \\
= & \
\sum_{k \in B_n(0)}
\left |
\alpha_k^{(\varepsilon)} e_k(x)
\right |^2
\frac{1}{2 \mu_k^2} \left(1- e^{-2t \mu_k^2} \right)
\leq 
C
\sum_{k \in B_n(0)}
\frac{1}{\mu_k^2}
< 
\infty 
\end{align*}
for each $t \in [0,T]$, $x \in [0,L]^2$, and $n \in \N$, where $C>0$ is a constant, as the sequence~$( \alpha_k^{(\varepsilon)} )_{k \in \Z}$ is uniformly bounded.
Similarly, for $n>m$, we can bound the difference of two elements of this sequence
\[  \mathbb{E} \left [ \left | Z^{(n)}(t,x)-Z^{(m)}(t,x) \right|^2 \right ]
\leq C \sum_{k\in B_n(0)\setminus B_m(0)}\frac1{\mu_k^2} 
\] 
 and thus show that 
$
\left(
Z^{(n)}(t,x)
\right)_{n \in \N}
$ is a Cauchy sequence in~$L^2(\Omega,\mathbb{R})$.
\end{proof}

Our next aim is to provide bounds on $Z$, $Z_\varepsilon$, $\nabla Z_\varepsilon$ and $Z-Z_\varepsilon$. Since all of these estimates are very similar, we first establish a more general result covering all these cases.

\begin{lemma}
\label{lem:techOU}
Let $\delta\in(0,1/4)$ and let
$\left(X_k\right)_{k\in\mathcal Z}$ be a family of independent,
centred, real-valued Gaussian processes on $[0,T]$.
Assume that there exist constants $C_0,C_\delta>0$ such that
\begin{align} \label{X:bounds}
\sup_{t\in[0,T]}
\mathbb E\left[|X_k(t)|^2\right]
&\leq
C_0\mu_k^{-2}, \qquad
\mathbb E\left[
|X_k(t)-X_k(s)|^2
\right]
 \leq
C_\delta
|t-s|^\delta
\mu_k^{4\delta-2} 
\end{align}
for all $k\in\mathcal Z$ and $s,t\in[0,T]$.
Let $(\alpha_k)_{k\in\mathcal Z}$ be a deterministic real-valued
sequence satisfying
\begin{align*} 
\sum_{k\in\mathcal Z}
\alpha_k^2\mu_k^{4\delta-2}
<\infty
\end{align*}
and define
\[
X_N(t,x)
\coloneqq
\sum_{k\in B_N(0)}
\alpha_k X_k(t)e_k(x), \qquad N \in \N,
\]
where $B_N(0) \subset \ZZ$ is the centered ball with radius $N \in \N$.
Then, for every $p\geq1$, the sequence $(X_N)_{N\in\mathbb N}$
converges to a process $X$ in
$
L^p\left(
\Omega ,
C^0\left([0,T]\times[0,L]^2\right)
\right)$.
In particular,
\begin{align} \label{gener:X}
    X(t,x)
=
\sum_{k\in\mathcal Z}
\alpha_kX_k(t)e_k(x)
\end{align} 
admits a continuous version and we have
\begin{align} 
\label{e:boundX}
\mathbb E
\left[
\|X\|_{C^0([0,T]\times[0,L]^2)}^p
\right]
\leq
C 
\left(
\sum_{k\in\mathcal Z}
\alpha_k ^2\mu_k^{4\delta-2}
\right)^{p/2}
\end{align}
where $C >0$ depends on $C_\delta$, $C_0$,
  $T$, $\delta$, $L$, and $p$.
\end{lemma}

\begin{proof}
 It suffices to prove the assertion for $p>6/\delta$, since the
remaining cases follow from Hölder's inequality. Fix therefore
$p>6/\delta$ 
and choose $\alpha \in \big( \frac{3}{p} , \frac{\delta}{2}\big)$.
Let $A\subset\mathcal Z$ be finite. In the following, $C>0$ denotes a
generic constant which may vary from line to line and may depend on $C_\delta$, $C_0$,
$\alpha$, $T$, $\delta$, $L$, and $p$.
For $t,s\in[0,T]$ and $x,y\in[0,L]^2$, Gaussian moment estimates
and the independence of the processes $(X_k)_{k\in\mathcal Z}$ yield
\[
\mathbb E 
\Big|
\sum_{k\in A}\alpha_k
\bigl(X_k(t)e_k(x)-X_k(s)e_k(y)\bigr)
\Big|^p 
 \leq
C_p\left(
\mathbb E
\Big|
\sum_{k\in A}\alpha_k
\bigl(X_k(t)e_k(x)-X_k(s)e_k(y)\bigr)
\Big|^2
\right)^{p/2}.
\]
Using
\[
X_k(t)e_k(x)-X_k(s)e_k(y)
=
\bigl(X_k(t)-X_k(s)\bigr)e_k(x)
+
X_k(s)\bigl(e_k(x)-e_k(y)\bigr),
\]
we obtain
\begin{align*}
&\mathbb E
\Big|
\sum_{k\in A}\alpha_k
\bigl(X_k(t)e_k(x)-X_k(s)e_k(y)\bigr)
\Big|^2
\\
& \quad \leq
C 
\sum_{k\in A}\alpha_k^2
\mathbb E\left[
|X_k(t)-X_k(s)|^2
\right]
+
C
\sum_{k\in A}\alpha_k^2
\mathbb E\left[
|X_k(s)|^2
\right]
|e_k(x)-e_k(y)|^2.
\end{align*}
Using the explicit form of the Fourier basis functions, we have
\[
|e_k(x)-e_k(y)|^2
\leq C_\delta \mu_k^\delta |x-y|^\delta.
\]
Therefore, using \eqref{X:bounds}, we derive
\begin{align*}
&\mathbb E
\Big|
\sum_{k\in A}\alpha_k
\bigl(X_k(t)e_k(x)-X_k(s)e_k(y)\bigr)
\Big|^2
\leq
C
\left(
|t-s|^\delta
\sum_{k\in A}\alpha_k^2\mu_k^{4\delta-2}
+
|x-y|^\delta
\sum_{k\in A}\alpha_k^2\mu_k^{\delta-2}
\right)
\\
&\quad\leq
C
\left(
|t-s|+|x-y|
\right)^\delta
\sum_{k\in A}\alpha_k^2\mu_k^{4\delta-2}.
\end{align*}
Consequently, we have
\begin{align*}
&\mathbb E
\Big|
\sum_{k\in A}\alpha_k
\bigl(X_k(t)e_k(x)-X_k(s)e_k(y)\bigr)
\Big|^p
 \leq
C
\left(
|t-s|+|x-y|
\right)^{\frac{\delta p}{2}}
\Big(
\sum_{k\in A}\alpha_k^2\mu_k^{4\delta-2}
\Big)^{\frac{p}{2}}.
\end{align*}
Similarly, we obtain
\begin{align*}
\mathbb E
\Big|
\sum_{k\in A}\alpha_kX_k(t)e_k(x)
\Big|^p
&\leq
C
\Big(
\sum_{k\in A}\alpha_k^2\mu_k^{-2}
\Big)^{p/2}
 \leq
C
\Big(
\sum_{k\in A}\alpha_k^2\mu_k^{4\delta-2}
\Big)^{p/2}.
\end{align*}
Since $\alpha p>3$, Morrey's inequality yields
\begin{align*}
 \mathbb E
\Big\|
\sum_{k\in A}\alpha_kX_k e_k
\Big\|_{C^0([0,T]\times[0,L]^2)}^p
& \leq
C 
\mathbb E
\Big\|
\sum_{k\in A}\alpha_kX_k e_k
\Big\|_{W^{\alpha,p}([0,T]\times[0,L]^2)}^p
\\
& \leq
C 
\Big(
\sum_{k\in A}\alpha_k^2\mu_k^{4\delta-2}
\Big)^{p/2},
\end{align*}
where the fractional Sobolev seminorm is finite since
$\alpha<\delta/2$. Indeed, the corresponding singular integral is
bounded by
\[
\int_{[0,T]\times[0,L]^2}
\int_{[0,T]\times[0,L]^2}
\left(
|t-s|+|x-y|
\right)^{p(\delta/2-\alpha)-3}
\,\dd(t,x)\,\dd(s,y)
<\infty.
\]
\noindent
Finally, Hölder's inequality implies the claim 
for each $p \geq 1$. 

Taking $A=B_N(0)\setminus B_M(0)$, we obtain
\begin{align*}
\mathbb E
\|X_N-X_M\|_{C^0([0,T]\times[0,L]^2)}^p
\leq
C 
\Big(
\sum_{k\in B_N(0)\setminus B_M(0)}
\alpha_k^2\mu_k^{4\delta-2}
\Big)^{p/2}.
\end{align*}
 Therefore, by the summability assumption, the sequence $(X_N)_{N\in\mathbb N}$ is a Cauchy sequence in
$
L^p\left(
\Omega , 
C^0([0,T]\times[0,L]^2)
\right)$.
Hence, it converges to a process $X$ in this space. By the definition
of $X_N$, the limit is represented by \eqref{gener:X}.
Finally, taking $A = B_N(0)$ in the estimate above and letting $N\to\infty$ yields the bound \eqref{e:boundX}.
\end{proof}

 Using Lemma \ref{lem:techOU}, we obtain the following result:

\begin{theorem}
The processes $Z$ and $Z_\varepsilon$ belong to $L^p \left( \Omega , C^0\left( [0,T]\times [0,L]^2 \right) \right)$ for each $p\geq1$ and 
for every small $\delta>0$ there is a constant depending also on $p$, $T$, and $L$ 
such that 
\[
\mathbb{E} \left\| Z_\varepsilon
    \right \|_{C^0}^p
    \leq
    C \left( \ 
    \sum_{k \in \Z} \left (\alpha_k^{(\varepsilon)} \right )^2
     \mu_k^{4\delta - 2} 
    \right)^{p/2}.
\]
\end{theorem}

\begin{proof}
As the $Z_k$ are independent centred Gaussian processes with
\[
\mathbb{E}\left[Z_k(t)^2\right]
=
\frac{1-e^{-2t\mu_k^2}}{2\mu_k^2}
\leq
\frac{1}{2\mu_k^2},
\]
it remains to verify the increment estimate required in
\Cref{lem:techOU}. 
Without loss of generality, let $t\geq s$. Then we have
\begin{align*}
Z_k(t)-Z_k(s)
&=
\int_s^t e^{-(t-\tau)\mu_k^2}\,\dd\beta_k(\tau)
 +
\left(e^{-(t-s)\mu_k^2}-1\right)
\int_0^s e^{-(s-\tau)\mu_k^2}\,\dd\beta_k(\tau).
\end{align*}
By the It\^o isometry and the independence of the Brownian increments,
\begin{align*}
\mathbb{E}\left[
|Z_k(t)-Z_k(s)|^2
\right] 
& =
\frac{1-e^{-2(t-s)\mu_k^2}}{2\mu_k^2}
+
\left(1-e^{-(t-s)\mu_k^2}\right)^2
\frac{1-e^{-2s\mu_k^2}}{2\mu_k^2}
\\
&\leq
C\mu_k^{-2}
\min\left\{1,|t-s|\mu_k^2\right\}.
\end{align*}
Hence, for every $\delta\in(0,1/4)$, we obtain
\begin{align*}
\mathbb{E}\left[
|Z_k(t)-Z_k(s)|^2
\right] 
&\leq
C_\delta
|t-s|^\delta
\mu_k^{2\delta-2}
 \leq
C_{\delta,L}
|t-s|^\delta
\mu_k^{4\delta-2}.
\end{align*}
Therefore, all assumptions of \Cref{lem:techOU} are satisfied.
For $Z$ we choose $\alpha_k=1$, while for $Z_\varepsilon$ we choose
$\alpha_k=\alpha_k^{(\varepsilon)}$.
Since $\delta<1/4$, we have
\[
\sum_{k\in\mathcal Z}\mu_k^{4\delta-2}<\infty,
\]
and the assertion follows from \Cref{lem:techOU}.
\end{proof}

Moreover, by studying the partial derivatives $D_i Z_\varepsilon$ in the same way, we obtain the following result.
\begin{theorem}
\label{thm:regnabZ}
For every $p\geq1$, the process $\nabla Z_\varepsilon$ belongs to $L^p \left( \Omega , C^0\left( [0,T]\times [0,L]^2 \right) \right)$ and 
for every~$\delta \in \big(0 ,  \frac{1}{4} \wedge \frac{\delta(\varepsilon)}{4}  \big)$ there is a constant $C> 0$ depending also on $p$, $T$, and $L$ 
such that 
\[
\mathbb{E} \left\| \nabla Z_\varepsilon
    \right \|_{C^0}^p
    \leq
    C \left( \ 
    \sum_{k \in \Z} \left(\alpha_k^{(\varepsilon)} \right)^2
     \mu_k^{4\delta - 1} 
    \right)^{p/2}.
\]
\end{theorem}

Finally, we study the difference 
\[
X_\varepsilon \coloneqq Z - Z_\varepsilon = \sum_{k \in \mathcal{Z} } \left(1-\alpha_k^{(\varepsilon)} \right )  Z_k(t) \ e_k(x).
\]
As before, Lemma \ref{lem:techOU} yields, for any small $\delta>0$,
\[
\mathbb{E} \left\| Z-Z_\varepsilon
    \right \|_{C^0}^p
    \leq
    C \left( \ 
    \sum_{k \in \Z} \left(1-\alpha_k^{(\varepsilon)} \right )^2
     \mu_k^{4\delta - 2} 
    \right)^{p/2}.
\]
 
Since the $\alpha_k^{(\varepsilon)}$ are uniformly bounded, dominated convergence yields:

\begin{theorem} \label{conv:Z}
 For all $T>0$ and $p \geq 1$, we have $Z_\varepsilon\to Z $
 in $L^p \left( \Omega , C^0\left( [0,T]\times [0,L]^2 \right) \right)$
 for $\varepsilon\to 0$.
\end{theorem}

\section{Existence and uniqueness of the mild solution}
\label{sec:4}
	
	\noindent
	The following result follows directly from \emph{Banach's fixed-point theorem}:
	
	\begin{theorem} \label {thm:exuniq}
		For $u_0\in\mathcal{H}^1$, the stochastic partial differential equation \eqref{e:SPDEreg} according to \Cref{def:Zeps} has a unique mild solution $u_\varepsilon$ in the sense of \Cref{def:mild} in the space~$\mathcal{H}^1$.
	\end{theorem}

	\begin{remark}
		It is possible to consider initial conditions with
		$u_0 \notin \mathcal{H}^1$ by employing weighted spaces in time,
		e.g.,
		\begin{align*}
			X_{\alpha , \delta} (0,T) =
			\left \{ u \in C^0 \left( (0,T ] , \mathcal{H}^\alpha \right) :
			\sup_{t \in (0,T]} t^\alpha \| u (t) \|_{\mathcal{H}^\alpha}
			< \infty
			\right \}.
		\end{align*}
		However, we refrain from giving further details here.
		For an application of such initial conditions, see \cite{BlRo:10}.
	\end{remark}
	
	\noindent
	We now prove \Cref{thm:exuniq} based on Banach's fixed-point theorem.
	First, we define the mapping 
	\begin{align*}
		\mathfrak{G} \left( u_{\varepsilon} \right) (t)
		\coloneqq
		e^{tA} u _0 + 
		\int_{0}^t
		e^{(t-s)A}   \mathfrak{F} \left(u_{\varepsilon} (s) \right) \dd s
		+ Z_\varepsilon(t).
	\end{align*}
	To apply Banach's fixed-point theorem, we first establish the following result.
	\begin{lemma}
		The mapping~$f \circ \nabla : \mathcal{H}^{1} \left( \T^2 \right)
		\to L^2 \left ( \T^2 , \R^2 \right)$, 
		$
		h \mapsto \frac{ \nabla h }{ 1 + \left | \nabla h \right |^2}
		$
		(cf.\ \Cref{def:mild}),
		is bounded and globally Lipschitz continuous with the Lipschitz constant $1$.
	\end{lemma}
	\begin{proof}
		The boundedness follows immediately from the boundedness of the function~$f : \R^2 \to \R^2$.
		To obtain the global Lipschitz continuity of $f$, we use its continuous differentiability. Therefore,
		by the \emph{multidimensional mean value theorem}, we obtain for each~$x,y \in \R^2$ that
		\begin{align*}
			\left |
			f(x) - f(y)
			\right |
			\leq
			\sup_{t \in [0,1]}
			\left \{
			\left \|
			D{f} \left(x + t (y-x) \right)
			\right \|
			\right \}
			|
			x - y
			|.
		\end{align*}
		where
		\begin{align*}
			D{f} (z) \coloneqq 
			\frac{1}{\left(1 + |z|^2 \right)^2}
			\left(
			\begin{array}{rr}
				1- z_1^2 + z_2 ^2 & -2 z_1 z_2 \\
				-2 z_2 z_1 & 1 + z_1^2 - z_2 ^2
			\end{array}
			\right) \in \R^{2 \times 2}
		\end{align*}
		denotes the Jacobian of the function~$f$ for $z \in \R^2$.
		We have $
			\left \|
			D {f}
			\right \|_{\infty} \leq 1$, 
		which implies the global Lipschitz continuity of~$f$.
		Applying Poincaré's inequality and the characterization of the Sobolev spaces $\mathcal{H}^{\alpha} = \mathcal{H}^{\alpha} \left(\T^2 \right)$ via the operator $(- \Delta)^{\frac{\alpha}{2}}$ (see \Cref{sec:2}),
		the linear operator $\nabla : \mathcal{H}^1 (\T^2, \R ) \to L^2 (\T^2, \R^2 )$ is a bounded linear operator with an operator norm $1$.
		Consequently, the mapping 
		\begin{align*}
			f \circ \nabla : \mathcal{H}^1 (\T^2, \R ) \to L^2 (\T^2, \R^2 ), \quad 
			h \mapsto f(\nabla h)
		\end{align*}
		is globally Lipschitz continuous.
		For each~$u,v \in \mathcal{H}^1$ we obtain the bound
		\begin{align*}
			\left \|
			\frac{\nabla u}{1 + \left | \nabla u \right |^2}
			 - 
			\frac{\nabla v}{1 + \left | \nabla v \right |^2}
			\right \|_{L^2}
			& \leq 
			\| u   -   v\|_{\mathcal{H}^1} .
		\end{align*}
	This verifies the lemma.
	\end{proof}

	\noindent
	Additionally, in preparation for applying Banach's fixed-point theorem, we recall several key properties of the semigroup $\bigl(e^{tA}\bigr)_{t \geq 0}$:
	
	\begin{itemize}
		\item \textbf{Smoothing estimate:}  
		For $\alpha < \beta$ and $t \in (0,T]$, the operator norm satisfies
		\begin{align} \label{Ineq:SG}
			\bigl\| e^{tA} \bigr\|_{L\!\left( \mathcal{H}^{\alpha}, \mathcal{H}^{\beta} \right)} 
			\leq 
			\left(\frac{\beta - \alpha}{4e}\right)^{ \tfrac{\beta - \alpha}{4} }
			t^{\tfrac{\alpha - \beta}{4}}.
		\end{align}
		This inequality is integrable at $t=0$ whenever $\alpha < \beta < \alpha+4$ (cf.~\eqref{e:SG}).
		
		\item \textbf{Divergence operator:}  
		The divergence operator
		\[
		\nabla \cdot : \mathcal{H}^1(\T^2, \R^2) \to L^2(\T^2,\R)
		\]
        is a bounded linear operator with an operator norm $1$.
		
		\item \textbf{Continuity of the semigroup:}  
		The operator $A$ generates an analytic semigroup~$\bigl(e^{tA}\bigr)_{t\geq 0}$ on $L^p(\T^2)$ for each $p \in [1,\infty)$.
        Moreover, the semigroup is strongly continuous on $\mathcal{H}^1$, and hence 
		\begin{align} \label{continofe}
			t \mapsto e^{tA} u_0 \in C^0([0,T], \mathcal{H}^{1}  )
		\end{align}
		for all $u_0 \in \mathcal{H}^{1}$ (see \cite[Section 2.6.]{PAZY}).
		
	\end{itemize}

	\begin{lemma}  \label{lem:continint}
		Let $w \in C^0\left([0,T], \mathcal{H}^1\right)$. Then the mapping
	\begin{align*}
		t
		\mapsto
		\int_{0}^t e^{(t-s)A} \mathfrak{F} \left(w(s) \right) \dd s
	\end{align*}
	belongs to $C^0 \left([0,T], \mathcal{H}^{3-} \right)$.
	\end{lemma}

	\begin{proof}
		Because the nonlinearity is Lipschitz continuous, $\mathfrak{F} : \mathcal{H}^1 \to \mathcal{H}^{-1}$ has Lipschitz constant $1$. Therefore, $\mathfrak{F}(w)$ belongs to $C^0([0,T], \mathcal{H}^{-1})$.
		By choosing $\alpha = -1$ and $\beta<3$ and applying \eqref{Ineq:SG},
		we obtain  
		the regularity
		\begin{align*}
			t \mapsto
			\int_{0}^t e^{(t-s)A} \mathfrak{F} \left(w(s) \right) \dd s
			\ \in    C^0 \left([0,T], \mathcal{H}^{3-} \right).
		\end{align*}
	The proof is now complete.
	\end{proof}
	
	\noindent
	Finally, by \eqref{continofe} and \Cref{lem:continint} we conclude that the operator~$\mathfrak{G}$ is a self-mapping on
	$C^0 \left([0,T], {\mathcal{H}^1}\right)$ provided that $Z_\varepsilon \in C^0 \left([0,T], {\mathcal{H}^1}\right)$, which we already established in \Cref{thm:regnabZ}.

	\begin{theorem}
		For sufficiently large $n\in\mathbb N$, the $n$-th iteration $\mathfrak{G}^{n} \coloneqq \overbrace{\mathfrak{G} \circ \cdots \circ \mathfrak{G}}^{n \ \text{times}}$ of the operator $\mathfrak{G}$ is a contraction on $C^0 \left([0,T], {\mathcal{H}^1} \right)$ and thus there is a pathwise unique fixed point.
		
	\end{theorem}
	
	\begin{proof}
	We consider the difference between two solutions $u, v \in C^0([0,T], {\mathcal{H}^1})$ with the same initial condition $u_0$.
	Thus, using \eqref{Ineq:SG} and the Lipschitz continuity of~$\mathfrak{F}$, we obtain
	\begin{align*}
		\left \|
		\mathfrak{G} (u)(t)
		-
		\mathfrak{G} (v)(t)
		\right \|_{\mathcal{H}^1}
		&=
		\left \|
		\int_{0}^t e^{(t-s)A} \mathfrak{F} \left(u(s) \right) \dd s
		-
		\int_{0}^t e^{(t-s)A} \mathfrak{F} \left(v(s) \right) \dd s
		\right \|_{\mathcal{H}^1}
		\\
		& \leq
		\sqrt{\frac{1}{2e}} \
		\int_{0}^t 
		\frac{1}{\sqrt{t-s} } 
		\left \|
		u(s) - v(s)\right \|_{\mathcal{H}^1} \dd s.
	\end{align*}
	By considering the supremum norm in time, we obtain
	\begin{align*}
		\sup_{t \in [0,T]}
		\left \{
		\left \|
		\mathfrak{G} (u)(t)
		-
		\mathfrak{G} (v)(t)
		\right \|_{\mathcal{H}^1}
		\right 
		\}
		& \leq 
		\sup_{t \in [0,T]}
		\left \{
		\sqrt{\frac{1}{2e}} 
		\int_{0}^t 
		\frac{1}{\sqrt{t-s} } 
		\left \|
		u(s) - v(s)\right \|_{\mathcal{H}^1} \dd s
		\right 
		\}
		\\
		& \leq
		\sqrt{\frac{1}{2e}} 
		\sup_{t \in [0,T]}
		\left \{
		\left \|
		u(t) - v(t)\right \|_{\mathcal{H}^1} 
		\right 
		\}
		\int_{0}^T 
		\frac{1}{\sqrt{s} }  
		\dd s
		\\
		& =
		\sqrt{\frac{2T}{e}} 
		\sup_{t \in [0,T]}
		\left \{
		\left \|
		u(t) - v(t)\right \|_{\mathcal{H}^1} 
		\right 
		\}.
	\end{align*}
	
	\noindent
	Since $\sqrt{2T/e}<1$ holds only for sufficiently small $T>0$, the
mapping $\mathfrak{G}$ is not necessarily a contraction on arbitrary
time intervals. We therefore consider its $n$-th iteration, where
$n\in\mathbb N$ is sufficiently large.
	For this purpose, we first compute the integral via substitution
	\begin{align*}
		\int_{0}^{t} \sqrt{ \frac{ z^{k} }{t - z}}  \dd z
		& =
		t^{\frac{k+1}{2}} \int_{0}^{1} \sqrt{ \frac{ x^{k} }{1 - x}}  \dd x 
		 =
		t^{\frac{k+1}{2}} \int_{0}^{1} x^{\frac{k}{2}} (1 - x)^{-\frac{1}{2}}  \dd x 
		 =
		t^{\frac{k+1}{2}}   B\left(\frac{k+2}{2} ,  \frac{1}{2} \right),
	\end{align*}
	for $k \in \N_0$,
	where
	$$
	B\left(
	y   ,   z \right) 
	\coloneqq 
	\frac{ \Gamma\left(	y \right) \Gamma\left( 	z \right)}{\Gamma\left(y+z\right)},
	$$
	for $y, z > 0$ is the~\emph{Beta function} and $\Gamma$ the \emph{Gamma function}.
	Furthermore, the resulting product is telescoping
	\begin{align*}
		\prod_{k=1}^{n-1}
		B\left(
		\frac{k+2}{2}   ,   \frac{1}{2} \right)
		=
		\prod_{k=1}^{n-1} 
		\frac{ \Gamma\left(	\frac{k+2}{2} \right) \Gamma\left( 	\frac{1}{2} \right)}{\Gamma\left(\frac{k+3}{2}\right)}
		= 
        \Gamma\left( 	\frac{1}{2} \right) ^{n-1}
		\frac{ \Gamma\left( 	\frac{3}{2} \right)}{ \Gamma\left( 	\frac{n+2}{2} \right)}.
	\end{align*}
	By combining these calculations, we can show that the $n$-th iteration satisfies the following contraction property
	\begin{align*}
		& \sup_{t \in [0,T]}
		\left \|
		\mathfrak{G}^n (u)(t)
		-
		\mathfrak{G}^n (v)(t)
		\right \|_{\mathcal{H}^1}
		\\
		& \leq 
		\left(
		\frac{1}{2e}
		\right)^{\frac{1}{2}}
        \sup_{t \in [0,T]}
		\int_{0}^t
		\frac{1}{\sqrt{t-s} } 
		\left \|
		\mathfrak{G}^{n-1} (u)(s)
		-
		\mathfrak{G}^{n-1}  (v)(s)
		\right \|_{\mathcal{H}^1}
		\dd s
		\\
		& \leq
		\left(
		\frac{1}{2e}
		\right)^{n/2} 
		\sup_{s \in [0,T]}
		\left \|
		u(s)
		-
		v(s)
		\right \|_{\mathcal{H}^1}
		\int_{0}^t
		\frac{1}{\sqrt{t-t_1} } 
		\hdots
		\int_{0}^{t_{n-1}}
		\frac{1}{\sqrt{t_{n-1}-t_n} } 
		\dd t_1 \hdots \dd t_n
		\\
		& \leq
		2
		\Gamma\left( 	\frac{3}{2} \right)
		\left(
		\frac{T}{2e}
		\right)^{n/2} 
        \Gamma\left( 	\frac{1}{2} \right)^{n-1} 
		\frac{1}{ \Gamma\left( 	\frac{n+2}{2} \right)} 
		\sup_{s \in [0,T]}
		\left \|
		u(s)
		-
		v(s)
		\right \|_{\mathcal{H}^1}.
	\end{align*}
		This expression converges to $0$ as $n \to \infty$.
		Thus, $\mathfrak{G}^n$ is a contraction for sufficiently large $n \in \mathbb{N}$.
		We now apply a corollary of Banach's fixed-point theorem \cite{WERN} for fixed-point iterations.
		Since the stochastic process $Z_\varepsilon$ is also an element of $C^0([0,T], {\mathcal{H}^1})$, we obtain the existence and uniqueness of a fixed point $\mathfrak{G}(u) = u$ in $C^0([0,T], {\mathcal{H}^1})$.
		\end{proof}
	
	\noindent
	By applying Banach's fixed-point theorem, we establish the existence and uniqueness of the mild solution $u_\varepsilon \in C^0([0,T], {\mathcal{H}^1})$, which for each $t \in [0,T]$ is given by
	\begin{align*}
		u_\varepsilon(t) =
		e^{tA} u _0 + 
		\int_{0}^t
		e^{(t-s)A} \mathfrak{F} \left(u_\varepsilon(s) \right)\dd s + Z_\varepsilon(t).
	\end{align*}

	\section{Uniform boundedness}
	 \label {sec:5}
	
	In the following, we establish a deterministic uniform bound on $\nabla v_\varepsilon$, which is crucial for proving the vanishing of the nonlinearity.
	For this purpose, consider the standard transformation 
	\[
	v_\varepsilon \coloneqq u_\varepsilon - Z_\varepsilon - e^{tA} u_0 ,
	\]
	where $v_\varepsilon$ solves the stochastic partial differential equation
	\begin{equation} \label {e:defv}
		\partial_t v_\varepsilon = - \Delta^2 v_\varepsilon - \nabla \cdot f 
		\left(
		\nabla  v _\varepsilon + \nabla Z _\varepsilon + \nabla e^{tA} u_0
		\right), \quad v_\varepsilon(0) = 0.
	\end{equation}
	
	\begin{theorem} \label{thm:boundv}
		For all $T>0$, 
		there is a constant $C = C_{L,T} >0$ such that 
		\begin{align*}
			|\nabla v_\varepsilon(\omega , t,x)| \leq C
		\end{align*}
		for almost every $\omega\in\Omega$, every $t\in[0,T]$, $x\in\T^2$, and $\varepsilon>0$.
	\end{theorem}

	\begin{proof}
		First, note that $\|f \|_{L^\infty \left(\R^2 , \R^2 \right)}$ is bounded by $1$.
		Let $\delta \in(0,1)$.
		For each $t \in [0,T]$ and $\varepsilon>0$, the function $v_\varepsilon$ satisfies
		\begin{align*}
			& \left \|\nabla v_\varepsilon(t)\right\|_{L^\infty \left( \T^2 , \R^2 \right)}
			\leq  C_\delta
			\left \|
			v_\varepsilon(t)
			\right \|_{\mathcal{H}^{2+\delta}}
			\\
			& = 
			C_\delta
			\left \|
			\int_{0}^t e^{(t-s)A}  \mathfrak{F} \left(v_\varepsilon(s) + Z_\varepsilon(s) + e^{sA} u_0 \right)   \dd s
			\right \|_{\mathcal{H}^{2+\delta}}
			\\
			& \leq
			C_\delta
			\int_{0}^t \left \|e^{(t-s)A} \nabla \cdot \right \|_{L\left( L^2\left(\T^2 , \R^2 \right) , \mathcal{H}^{2+\delta} \right)} \left \| f \left(\nabla v_\varepsilon(s) + \nabla Z_\varepsilon(s) + \nabla e^{sA} u_0 \right)
			\right \|_{L^2} \dd s
			\\ & \leq
			C_\delta\int_{0}^t  
			(t-s)^{-(3+\delta) /4}   L   \dd s
			\leq 
			C_{\delta}
			 L
			 T^{(1-\delta)/4}.
		\end{align*} 
	This establishes the desired result.
	\end{proof}
	
	\section{Vanishing of the nonlinearity}
	 \label {sec:6}
	
	In this subsection, we prove the convergence $f \left(\nabla u_\varepsilon(t,x) \right) \to 0$ for $\varepsilon \to 0$ in $L^p\left(\Omega, \R^2\right)$, uniformly in $x \in \T^2$ for all $t \in (0,T]$.
	First, we present a short auxiliary statement that relies on the periodicity of the domain and the type of regularization.
	\begin{lemma}  \label {lem:indepx}
		For each~$t \in [0,T]$ the covariance operators of the processes~$Z_\varepsilon(t,x)$ and $\nabla Z_\varepsilon(t,x)$ do not depend on~$x \in \T^2$.
	\end{lemma}
	
	\begin{proof}
		Without loss of generality, assume $k_1 , k_2 > 0$.
     Applying the \emph{Pythagorean identity} to the definition of $e_k$ in \Cref{sec:2}, we obtain for each $x = \left(x_1 , x_2 \right) \in \T^2$
		\begin{align*}
			& e_{\left(k_1 , k_2 \right)}^2 \left( x \right) 
			 +
			e_{\left(- k_1 , k_2 \right)}^2 \left( x \right) 
			+
			e_{\left(k_1 , -k_2 \right)}^2 \left( x \right) 
			+
			e_{\left(- k_1 , - k_2 \right)}^2 \left( x \right) 
			\\
			& = 
			\frac{4}{L^2} \left[
			\sin^2 \left( \frac{2 \pi k_1}{L} x_1 \right) \sin^2 \left( \frac{2 \pi k_2}{L} x_2 \right)
			+
			\sin^2 \left( \frac{2 \pi k_1}{L} x_1 \right) \cos^2 \left( \frac{2 \pi k_2}{L} x_2 \right) \right]
			\\ & \quad +
			\frac{4}{L^2} \left[ 
			\cos^2 \left( \frac{2 \pi k_1}{L} x_1 \right) \sin^2 \left( \frac{2 \pi k_2}{L} x_2 \right)
			+
			\cos^2 \left( \frac{2 \pi k_1}{L} x_1 \right) \cos^2 \left( \frac{2 \pi k_2}{L} x_2 \right)
			\right]
			\\
			& = \frac{4}{L^2}
		\end{align*}
		which does not depend on $x \in \T^2$.
		Since $\omega_0(x_i) = L^{- \frac{1}{2}}$ is a constant function, the other cases with $k_1 = 0$ or $k_2 = 0$ follow analogously. Note that 
        \[
			\cos^2 \left( \frac{2 \pi l }{L} x_1 \right) + \sin^2 \left( \frac{2 \pi l}{L} x_1 \right)
			+
			\cos^2 \left( \frac{2 \pi l}{L} x_2 \right) + \sin^2 \left( \frac{2 \pi l}{L} x_2 \right) = 2.
        \]
		By the independence of the sequence of Wiener processes $\left(\beta_k\right)_{k \in \ZZ}$, we can apply the Itô isometry. Furthermore, due to the radial symmetry of the sequences $\left(\mu_k\right)_{k \in \ZZ}$ and $\big(\alpha_k^{(\varepsilon)}\big)_{k \in \ZZ}$, we conclude for each $\varepsilon > 0$ that
		\begin{align*}
			& \mathbb{E}
			\left[ 
			\left |
			Z_\varepsilon (t,x)
			\right |^2 
			\right ]   =
			\mathbb{E}
			\left[ 
			\left(
			\sum_{k \in \ZZ}
			 \alpha_k^{(\varepsilon)}  e_k(x) 
			\int_0^t e^{-(t-s) \mu_k^2} \dd \beta_k(s)
			\right)^2
			\right] 
			\\
			& =
			\sum_{k \in \ZZ}
			\left( \alpha_k^{(\varepsilon)} \right)^2  e_k^2 (x) 
			\int_0^t e^{-2(t-s) \mu_k^2} \dd s
			\\
			& = 
			\frac{4}{L^2}
			\sum_{\substack{k \in \ZZ \\ k_1, k_2 > 0}}
			\left( \alpha_k^{(\varepsilon)} \right)^2 
			\int_0^t e^{-2(t-s) \mu_k^2} \dd s
		  + 
			 \frac{4}{L^2}
			\sum_{l \in \N}
			\left( \alpha_{(l,0)}^{(\varepsilon)} \right)^2 
			\int_0^t e^{-2(t-s) \mu_{(l,0)}^2} \dd s,
		\end{align*}
		which does not depend on~$x \in \T^2$.
        Furthermore, 
        it is straightforward to verify, by applying the same argument to the derivatives of the Fourier basis functions, that the covariance operator of $\nabla Z_\varepsilon (t,x)$ is also independent of $x \in \T^2$.
        This establishes the desired result.
	\end{proof}

	\begin{theorem}  \label {UCV} 
		For each $p \geq 1$ and $t \in (0,T]$, we have
			\begin{align*}
			\sup_{x \in \T^2 } 
			\mathbb{E} \left[ \left | 
			f \left(\nabla u_\varepsilon (t,x) \right) 
			\right |^p \right] \to 0 \quad \text{if }  \varepsilon \to 0.
		\end{align*}
	\end{theorem} 

\begin{remark}
    Although the analyticity of $\left(e^{tA} \right)_{t \geq 0}$ was stated above only on $L^p(\T^2)$ for $p \in (1, \infty)$, the semigroup is uniformly bounded on $W^{1 , \infty} \left(\T^2 \right)$ on bounded time intervals.
    In particular, for $u_0\in W^{1,\infty}(\T^2)$, we have
$$ \|  \nabla e^{\cdot  A} u_0 \|_{L^\infty \left([0,T] \times \T^2 \right)}  < \infty. $$
\end{remark}
    
	\begin{proof}
		Fix an arbitrary $t \in (0,T]$.
		Due to  \cref{thm:boundv} we fix a constant $M>0$ such that 
		$M\geq \| \nabla v_\varepsilon\|_{L^\infty \left([0,T] \times \T^2 \right)} + \|  \nabla e^{\cdot  A} u_0 \|_{L^\infty \left([0,T] \times \T^2 \right)} $.
		We also choose $K_\varepsilon \geq M$ by
		\begin{align*}
			K_\varepsilon \coloneqq
			\left[
			\sum_{k_1 , - k_2 \in \mathbb{N}
			} \left( \frac{ k_1  }{ | k|^2  }\right)^2
			\left(	\alpha_k^{(\varepsilon)} \right)^2
			\right]^\eta \vee M,
		\end{align*}
		where $\eta$ is a constant in the interval $\left(0 , \frac{1}{2} \right)$. Thus, we derive
		$
		K_\varepsilon \to \infty
		$
		as $\varepsilon \to 0$.
		First, note that we can directly estimate
		\[ |f(\nabla u_\varepsilon (t,x)) |
		\leq \frac{2}{1+|\nabla u_\varepsilon (t,x)|}.
		\]
		As $\nabla u_\varepsilon=\nabla v_\varepsilon + \nabla Z_\varepsilon + \nabla e^{t A} u_0$, we obtain 
		for $| \nabla Z_\varepsilon|>K_\varepsilon$ the inequality
		\[
		\frac{2}{1 + \left | \nabla u_\varepsilon (t,x)\right | } 
		\leq   
		\frac{2}{1 + \left |  \nabla Z_\varepsilon (t,x)\right | - \left |  \nabla v_\varepsilon (t,x)\right | 
        - \left |  \nabla e^{t A} u_0 (x) \right | 
        }
		\leq  \
		\frac{2}{1 + K_\varepsilon - M}
		\]

		\noindent
		Let $B_1(0)$ be the ball with radius $1$ in $\R^2$.
		By the boundedness of $f$, we obtain 
		\begin{align*}
			\mathbb{E} \left[ \left | 
			f \left(\nabla u_\varepsilon (t,x)\right) 
			\right |^p \right]
			& \leq 
			\mathbb{E} \left[ \left | 
			f \left(\nabla u_\varepsilon (t,x)\right) 
			\right |^p
			\mid  | \nabla Z_\varepsilon (t,x) | > K_\varepsilon
			\right]   \mathbb{P} \left( | \nabla Z_\varepsilon(t,x) | > K_\varepsilon \right)
			\\
			&   \  +
			\mathbb{E} \left[ \left | 
			f \left(\nabla u_\varepsilon (t,x)\right) 
			\right |^p
			\mid  | \nabla Z_\varepsilon (t,x) | \leq K_\varepsilon
			\right]   \mathbb{P} \left( | \nabla Z_\varepsilon (t,x)| \leq K_\varepsilon \right)
			\\
			&
			\leq
			\left( \frac{2}{1 + K_\varepsilon - M} \right)^p + C   \mathbb{P} \left( | \nabla Z_\varepsilon (t,x)| \leq K_\varepsilon \right)
			\\
			&
			=
			\left( \frac{2}{1 + K_\varepsilon - M} \right)^p 
			+ C   \mathbb{P} \left( \frac{ \nabla Z_\varepsilon (t,x)}{K_\varepsilon} \in \overline{B_1(0) } \right)
			\\
			&
			\leq
			\left( \frac{2}{1 + K_\varepsilon - M} \right)^p 
			+ C   
			\max_{y \in \overline{B_1(0) }} \{ \varphi_\varepsilon (t,x, y) \} 
			\operatorname{vol} \left( B_1(0)  \right),
		\end{align*}
		where  $C \in (0,1)$ is a constant and
		$$ \varphi_\varepsilon(t,x,y)
		\coloneqq 
		\frac{1}{ 2 \pi \sqrt{\operatorname{det}(\Sigma_\varepsilon(t,x))}} 
		e^{
		- \frac{1}{2}   y^T   \Sigma_\varepsilon^{-1}(t,x) y},
		$$
		is the density function of
		$
		\frac{ \nabla Z_\varepsilon(t,x)}{K_\varepsilon} \sim
		\mathcal{N} \left(\mathbf{0} , \Sigma _\varepsilon (t,x) \right)$.
		\noindent
		We denote the covariance matrix by $\Sigma_\varepsilon$ as follows:
		\begin{align*}
			\Sigma_\varepsilon (t,x) 
			& \coloneqq   \operatorname{Cov} \left(  \frac{ \nabla Z_\varepsilon(t,x) }{K_\varepsilon} \right)
			=
			\frac{1}{K^2_\varepsilon}
			\sum_{k \in \cZ} 
			\left[ 
			1 - e^{-2t \mu_k^2}
			\right ] 
			\frac{\left( \alpha_k^{(\varepsilon)}\right)^2}{2 \mu_k^2}  \nabla e_k(x) \nabla e_k(x)^T
			\\ & =
			\frac{1}{K^2_\varepsilon}
			\sum_{k \in \cZ} 
			\left[ 
			1 - e^{-2t \mu_k^2}
			\right ] 
			\frac{\left( \alpha_k^{(\varepsilon)}\right)^2}{2 \mu_k^2} \
			\left(
			\begin{array}{rr} D_1 e_k(x) D_1 e_k(x) & D_1 e_k(x) D_2 e_k(x) 
				 \\ 
				D_2 e_k(x) D_1 e_k(x) & D_2 e_k(x) D_2 e_k(x)
			\end{array} \right).
		\end{align*}
		By \Cref{lem:indepx}, the covariance of  $Z_\varepsilon$, and thus the covariance of~$\nabla Z_\varepsilon$, does not depend on the variable~$x \in \T^2$.
		We investigate the determinant of the covariance matrix at $(0,0) \eqqcolon \mathbf{0}$, which simplifies the argument.
		Furthermore, by the definition of~$\left( e_k \right)_{k \in \ZZ}$, we obtain
		\begin{align*}
			D_1 e_k \left(\mathbf{0} \right) D_2 e_k \left( \mathbf{0} \right) 
			=
			\frac{2}{L} \left(\frac{2 \pi}{L}\right)^2 k_1 k_2
			\sin (0)
			\cos (0)
			\sin (0)
			\cos (0)
			=
			0
		\end{align*}
		for each~$k = \left(k_1 , k_2 \right  ) \in \ZZ$ and 
		thus we derive
		\begin{align}  \label {det:cov}
			\lefteqn{
				\operatorname{det}\left (\Sigma_\varepsilon (t,\mathbf{0}) \right) } 
			 \\
			& =  
			\frac{
				1  }{ 4 K^4_\varepsilon}
			\sum_{k, l  \in \cZ}
			\left( \frac{ \alpha_k^{(\varepsilon)} \alpha_l^{(\varepsilon)} }{\mu_k    \mu_l}\right)^2 
			\left[ 
			1 - e^{-2t \mu_k^2}   \right ] \left[ 
			1 - e^{-2t \mu_l^2}   \right ] \nonumber 
			\\ & \quad \quad  \left[
			\left | D_1 e_k \left( \mathbf{0} \right) \right |^2 
			\left | D_2 e_l \left( \mathbf{0} \right) \right |^2 
			-
			D_1 e_k \left( \mathbf{0} \right) D_1 e_l \left( \mathbf{0} \right) D_2 e_k \left( \mathbf{0} \right) D_2 e_l \left( \mathbf{0} \right)
			\right]	 \nonumber 
			\\
			&	=
			\frac{1}{4 K^4_\varepsilon}
			\sum_{k, l  \in \cZ} \left( \frac{ \alpha_k^{(\varepsilon)} \alpha_l^{(\varepsilon)} }{\mu_k    \mu_l}\right)^2
			\left | D_1 e_k \left( \mathbf{0} \right) \right |^2 
			\left | D_2 e_l \left( \mathbf{0} \right) \right |^2  \left[ 1 - e^{-2t \mu_k^2} \right] \left[ 
			1 - e^{-2t \mu_l^2}   \right ] \nonumber 
			\\ & = 
			\frac{1}{4 K^4_\varepsilon}
			\left[
			\sum_{k  \in \cZ} \left( \frac{ \alpha_k^{(\varepsilon)} }{\mu_k  }\right)^2
			\left | D_1 e_k \left( \mathbf{0} \right) \right |^2 \left[ 1 - e^{-2t \mu_k^2} \right]
			\right]^2 \nonumber 
			\\ & =
			\frac{1}{4 K^4_\varepsilon}
			\left[
			\sum_{\substack{k_1 \in \mathbb{N} \\ -k_2 \in \mathbb{N}_0}} \left( \frac{ \alpha_k^{(\varepsilon)} }{\mu_k  }\right)^2
			\frac{8 \pi ^2}{L^3}    k_1^2     \omega_{k_2} (0)^2 \left[ 1 - e^{-2t \mu_k^2} \right]
			\right]^2 \nonumber 
			\\ & \asymp  
			\frac{ 4}{K^4_\varepsilon \pi^4}
			\left[
		\sum_{\substack{k_1 \in \mathbb{N} \\ -k_2 \in \mathbb{N}_0}}  \frac{   k_1^2 }{|k|^4 }
			\ \left( \alpha_k^{(\varepsilon)}\right) ^2  
			\left[ 1 - e^{-2t \mu_k^2} \right]
			\right]^2 , \nonumber
		\end{align} 
		as $\omega_{k_2} (0) \in \big \{ L^{-\frac{1}{2}}  \ ,\ \sqrt{2}  L^{-\frac{1}{2}}  \big \}$
        and $
		1 - e^{-2t \left(2 \pi / L \right)^4}
		\leq 
		1 - e^{-2t \mu_k^2}
		\leq 1$
        for each~$t \geq 0$ and $k = (k_1, k_2) \in \cZ$. 
		For $t = 0$, this determinant is equal to~$0$.
		For fixed $t > 0$, the convergence of this series is equivalent to the convergence of the following series due to the non-negativity of its summands
		\begin{align*}
			\frac{1}{K_\varepsilon^2} \
			\sum_{k_1 , - k_2 \in \mathbb{N}
			} \left( \frac{ k_1  }{ | k|^2  }\right)^2
			\left(	\alpha_k^{(\varepsilon)} \right)^2 
            \asymp  
            \left(
			\sum_{k_1 , - k_2 \in \mathbb{N}
			} \left( \frac{ k_1  }{ | k|^2  }\right)^2
			\left(	\alpha_k^{(\varepsilon)} \right)^2
            \right)^{1- 2 \eta}
		\end{align*}
		which is finite for every fixed $\varepsilon>0$ due to the regularization.
    As $\varepsilon\to0$, Fatou's lemma and
    $\alpha_k^{(\varepsilon)}\to1$ imply that the series diverges.
    Since $1-2\eta>0$, the expression above therefore diverges as well.
		This follows directly from the \emph{integral comparison criterion} by analysing the convergence or divergence of
		\begin{align*}
			\int_{1}^\infty \int_{1}^\infty \frac{x_1^2}{\left( |x_1| + |x_2| \right)^4}  \dd x _1 \dd x _2
			& = 	\int_{1}^\infty 
			\frac{ 3 + 3 x_2 + x_2^2 }{ 3 (1 + x_2)^3}   \dd x_2 
			= \infty.
		\end{align*}
		Since this series diverges, we conclude that
		\begin{align*}
			\mathbb{E} \left[ \left | 
			f \left(\nabla u_\varepsilon (t,x) \right) 
			\right |^p \right]
			& \leq 
			\left( \frac{2}{1 + K_\varepsilon - M } \right)^p 
			+ C
			\operatorname{vol} \left( B_1(0)  \right)
			\max_{y \in \overline{B_1(0) }} \{ \varphi_\varepsilon (t,x, y) \} 
			\\
			& \leq
			\left( \frac{2}{1 + K_\varepsilon - M} \right)^p 
			+  \frac{C}{\sqrt{\operatorname{det}\left(\Sigma_\varepsilon (t,x) \right)}}
			\xrightarrow{\varepsilon \to 0} 0
		\end{align*}
		also vanishes in the limit for $t \in (0,T]$ and  $x \in \T^2$.
	\end{proof}

	\noindent
	Therefore, by pointwise convergence and the uniform boundedness of $f$, we directly conclude the following corollary.
	\begin{corollary}  \label {LPC}
		For each $p \geq 1$, we have
		$$
		\mathbb{E} \left[ \left \| 
		f \left(\nabla u_\varepsilon  \right) 
		\right \|^p_{L^p \left( [0,T]\times \T^2 , \R^2\right) } \right] \to 0
		\quad\text{if } \varepsilon \to 0.$$
	\end{corollary}
	
	\begin{proof}
		By the pointwise convergence in \Cref{UCV}
		\begin{align*}
			\mathbb{E} \left[ \left | 
			f \left(\nabla u_\varepsilon (t,x) \right) 
			\right |^p \right] 
			\xrightarrow{\varepsilon \to 0} 0
		\end{align*}
		for each $x \in \T^2$ and $t \in (0,T]$ and  the uniform boundedness
		$$
		\left \|
		f \left(
		\nabla u_\varepsilon
		\right)
		\right \|_{L^\infty \left( 
			[0,T] \times \T^2
			\right)}
		\leq 1,
		$$ 
		\noindent
		we directly derive
		the following result  by applying Tonelli's theorem:
		\begin{align*}
			\mathbb{E} \left \| 
			f \left(\nabla u_\varepsilon  \right) 
			\right \|^p_{L^p \left( [0,T]\times \T^2 , \R^2\right) }
			 & 
			=
			\mathbb{E} 
			\int_{[0,L]^2} \int_0^T 
			\left |
			f \left(\nabla u_\varepsilon (t,x) \right) 
			\right |^p \dd t \dd x 
			\\ & 
			= C 
			\int_{[0,L]^2} \int_0^T 
			\mathbb{E} 
			\left |
			f \left(\nabla u_\varepsilon (t,x)  \right) 
			\right |^p 
			\dd t  \dd x 
				\\ & 
			\leq C
			\int_0^T
			\sup_{x \in \T^2 } 
			\mathbb{E} \left[ \left | 
			f \left(\nabla u_\varepsilon (t,x) \right) 
			\right |^p \right] 
			\dd t 
		\end{align*}
		where $C=C_{L,T}>0$ is a constant.
		This last expression tends to zero as $\varepsilon$ tends to zero using the dominated convergence theorem.
	 	This finishes the proof.
	\end{proof}

	\section{Proof of the main result}
	 \label{sec:7}
	Finally, we are able to prove the convergence of the mild solutions from \Cref{thm:conv} using \Cref{LPC}.
	Let us first consider $v_\varepsilon = u_\varepsilon - Z_\varepsilon - e^{tA} u_0$, as the mild solution of \eqref{e:defv}.
	In order to verify the main result of this chapter, we have to prove the following convergence for each $p \geq 1$:
	$$
	v_\varepsilon
	\xrightarrow{\varepsilon \to 0} 0
	\quad \text{
		in }\quad 
	L^p\left(\Omega, C^0( [0,T] \times \T^2 ) \right).$$
	
	\begin{proof}[Proof of \Cref{thm:conv}]
		For each $p>2$, there exists a $\delta \in (0,1)$ such that $p > \frac{4}{2-\delta}$ or, equivalently,
		$
		\frac{(2+\delta)p}{4(p-1)} \in(0,1).
		$
		Now, by applying Morrey's inequality (cf.\ \cite[p.283]{LEPDE}) and Hölder's inequality, we obtain
		\begin{align*}
			&
			\mathbb{E} \left[\left \| v_\varepsilon   \right \|^p_{C^0 \left([0,T]\times \T^2 \right)}
			\right]  
			 \leq 
			 C_{T,L, \delta} \mathbb{E}\left[
			 \sup_{t\in[0,T]} \left \| v_\varepsilon (t) \right \|^p_{\mathcal{H}^{1+\delta}} 
			 \right]
			\\ & \leq C_{T,L, \delta}   
			\mathbb{E} \left[ \sup_{t\in[0,T]} \left \| 
			\int_0^t
			e^{(t-s)A} \nabla \cdot f \left( \nabla u_\varepsilon (s) \right) \dd s
			\right \|^p_{\mathcal{H}^{1+\delta} } 
			\right]
			\\
			& \leq C_{T,L, \delta}   
			\mathbb{E} 
			\left[
			\sup_{t\in[0,T]} \left( 
			\int_0^t (t-s)^{\frac{-(2+\delta)}{4}}
			\left\| f \left( \nabla u_\varepsilon (s) \right) \right\|_{L^2} \dd s
			\right)^p 
			\right]
			\\
			& \leq C_{T,L, \delta}   
			\mathbb{E} \left[ \sup_{t\in[0,T]} \left( 
			\int_0^t (t-s)^{\frac{-(2+\delta)p}{4(p-1)}} \dd s
			\right)^{(p-1)}   
			\int_0^t 
			\left\| f \left( \nabla u_\varepsilon (s) \right) \right\|^p_{L^2} \dd s
			\right]
			\\
			& \leq C_{T,L, \delta}   
			\mathbb{E} 
			\left[
			\int_0^T 
			\left\| f \left( \nabla u_\varepsilon (s) \right) \right\|^p_{L^2} \dd s
			\right]
			\\
			&  \leq C_{p, T,L, \delta}   
			\mathbb{E} 
			\left[
			\left\| f \left( \nabla u_\varepsilon \right) \right\|^{p}_{L^p([0,T]\times \T^2)} 
			\right]
			\xrightarrow[\varepsilon \to 0]{} 0
		\end{align*}
		which follows immediately from \cref{LPC}.
		Finally, by applying Hölder's inequality, the convergence
\[
v_\varepsilon \xrightarrow{\varepsilon\to0} 0
\]
holds in $L^p(\Omega,C^0([0,T]\times\T^2))$ for each $p\geq1$.
Since $
u_\varepsilon-u
=
v_\varepsilon + Z_\varepsilon-Z
$,
the convergence $Z_\varepsilon\to Z$ from \Cref{conv:Z} yields
\[
u_\varepsilon\to u
\quad\text{in}\quad
L^p(\Omega,C^0([0,T]\times\T^2)).
\]
This establishes the desired result.
	\end{proof}
	
	\noindent
	This convergence yields the reduction to the linear limiting SPDE
\[
\partial_t u = - \Delta^2 u + \partial_t W.
\]
The linear instability of the original model and the corresponding
spinodal decomposition are investigated in more detail in
\cite{BR:GR26}.

\paragraph{Acknowledgements} 
Both authors acknowledge the support of the DFG through grant BL 535/12-1, project number 514726621.

\end{document}